\documentstyle[12pt]{article}
\title{\Large Some equivalent of extremally disconnected spaces}
\author{Vishvajit V S Gautam \thanks{Supported by UGC, New Delhi.\,\,\, }\\
{\sf Department of Mathematics}\\ 
{\sf  Kurukshetra University, Kurukshetra - 136 119 India.}}
\date{ }
\newtheorem{lemma}{\large \bf Lemma}
\newtheorem{coro}[lemma]{\large \bf Corollary}
\newtheorem{th}[lemma]{\large \bf Theorem}
\begin{document}
\maketitle
\begin{abstract}
In this short note we give some equivalent characterizations of extremally
disconnected spaces.\\\
\end{abstract}

Extremally disconnected spaces play a prominent role in set-theoretical
topology. Due to their peculiar properties extremally disconnected spaces
provide curcial applications in the theory of Boolean algebra, in axiomatic
set theory and in some branches of functional analysis (for example: in
$C^{*}$-algebra) as well. There are many interesting equivalent of  
extremally disconnected spaces. See for example : problems 1H, 3N, 6M of [3]
and, in case of Boolean algebra, Theorem 2.33 of [2]. \\
In this short note we prove some set-theoretical equivalent of extremally 
disconnected spaces.
Proofs of these equivalent conditions are not trivial. But as a good
exercies details are left to the readers.\\
Recall that a topological space $X$ is said to be extremally disconnected if
the closure of every open set of $X$ is open in $X$.\\
We denote by $A^{o}$ the interior of A and by $A^{-}$ the closure of A.

\begin{lemma}
{\em Let $X$ be a topological space. If $A$ is a subset of $X$ and $B$ is an
open subset of $X$ then
\[A^{-o} \cap B^{-o} = (A \cap B)^{-o}. \] 
}
\end{lemma}
{\large \bf Proof.}  Use the fact $A \cap B^{-o} \subset (A \cap B)^{-}$. $\Box$\\

As a Corollary to this we prove the Lemma 2.34 of [2].
\begin{coro}
{\em  If $A$ and $B$ are open subsets of a topological space $X$ then
\[A^{-o} \cap B^{-o} = (A \cap B)^{-o}. \] 
}
\end{coro}

\begin{th}
{\em  The following are equivalent for a space $X$ :\\
(a) $X$ is extremally disconnected;\\
(b) \(A^{-} \cap B^{-} = (A \cap B)^{-} \) for all open subsets
$A$ and $B$ of $X$;\\
(c) $A^{-} \cap B^{-} = \emptyset$ for all open subsets
$A$ and $B$ of $X$ with $A \cap B = \emptyset$;\\
(d) $K^{-o-} \cap A^{-} = \emptyset$ for all subsets
$K$ and all open subsets $A$ of $X$ with $K \cap A = \emptyset$;\\
(e) $K^{-o} \cap A^{-} = (K \cap A)^{-o}$ for all subsets
$K$ and all open subsets $A$ of $X$ with $K \cap A = \emptyset$;\\
(f) \(A^{-o} \cup B^{-o} = (A \cup B)^{-o} \) for all open subsets
$A$ and $B$ of $X$; and\\ 
(g) \(G^{o} \cup H^{o} = (G \cup H)^{o} \) for all closed subsets
$G$ and $H$ of $X$.
          
}
\end{th}
{\large \bf Proof.} (a) $\Longrightarrow$ (b) $\Longrightarrow$ (c)
$\Longrightarrow$ (a) $\Longleftrightarrow$ (e) $\Longrightarrow$ (d)
$\Longrightarrow$ (c) $\Longrightarrow$ (a) $\Longrightarrow$ (f)
$\Longrightarrow$ (g) $\Longrightarrow$ (c) $\Longrightarrow$ (a). $\Box$\\
\\
\
{\large \bf Remark.} Above Theorem seems to be an isolated result, but this
will be useful to the readers those are interested in generalization of concepts
having closed relation with extremally disconnected spaces. For example: in
generalization of certain properties  of spaces with filters and  $\infty$-disconnected (extremally
disconected) spaces (cf. [1]).

\section*{References.}
[1]  R.N. Ball, A.W. Hagner and A.J. Macula, An $\alpha$-disconnected space
     has no proper monic preimage, Topology Appl. 37(1990) 141-151.\\
\newline [2]    W.W. Comfort and S. Negrepontis, The Theory of Ultrafilters,
       Die Grundlehren der mathematischen Wissenschaften, Band 211
       (Springer-Verlag, NY, 1974).\\
\newline [3]  L. Gillman and M. Jerison, Rings of Continuous Functions
        (Van Nostrand Reinhold, Princeton, NJ, 1960).

\end{document}